\documentclass{article}
\usepackage{graphicx}
\usepackage{psfrag}
\usepackage{amssymb}
\usepackage{amsmath}
\usepackage{amsthm}
\usepackage{amscd}

\newtheorem{thm}{Theorem}[section]

\newtheorem{lem}[thm]{Lemma}

\newtheorem*{pnt}{Prime Number Theorem}

\theoremstyle{definition}
\newtheorem{defn}[thm]{Definition}

\title{On the Sum of Reciprocals of Primes}
\author{Young Deuk Kim
\\SNU College \\
Seoul National University\\
Seoul 08826, Korea
\\(yk274@snu.ac.kr)}
\date{\today}

\begin{document}
\maketitle

\begin{abstract}
Suppose that $y>0$, $0\leq\alpha<2\pi$, and $0<K<1$. Let $P^+$ be the set of primes $p$ such that $\cos(y\ln p+\alpha)>K$, and $P^-$ the set of primes $p$ such that $\cos(y\ln p+\alpha)<-K$. In this paper, we prove $\sum_{p\in P^+}\frac{1}{p}=\infty$ and
$\sum_{p\in P^-}\frac{1}{p}=\infty$.

\vspace{.3cm}
\end{abstract}

%%%%%%%%%%%%%%%%%%%%%%
%%%%%%%%%%%%%%%%%%%%%%
\section{Introduction}
%%%%%%%%%%%%%%%%%%%%%%
%%%%%%%%%%%%%%%%%%%%%%

Let $P$ be the set of primes and $\mathbb{N}$ be the set of natural numbers.
In 1737, Euler \cite{Euler} proved that the sum of reciprocals of primes diverges:
$$\sum_{p\in P}\frac{1}{p}=\infty.$$

\begin{defn}
Suppose that $y>0$, $0\leq\alpha<2\pi$, and $0<K<1$. Let
$$P^+(y,\alpha,K)=\{p\in P\mid \cos(y\ln p+\alpha)>K\}$$
and
$$P^-(y,\alpha,K)=\{p\in P\mid \cos(y\ln p+\alpha)<-K\}.$$
We write $P^+$ and $P^-$ for the sake of simplicity.
\end{defn}

\indent
Throughout this paper, we assume that $y>0$. We prove
\begin{thm}\label{thm}
$$\sum_{p\in P^+}\frac{1}{p}=\infty\quad\mbox{and}\quad \sum_{p\in P^-}\frac{1}{p}=\infty.$$
\end{thm}

\noindent
This theorem plays a crucial role in the author's works (\cite{kim},\cite{kim2}) related to the Riemann hypothesis.

%%%%%%%%%%%%%%%%%%%%%%%%%%%%%%%%%%%%%%%%%%%%%%%%%%%%
%%%%%%%%%%%%%%%%%%%%%%%%%%%%%%%%%%%%%%%%%%%%%%%%%%%%
\section{Proof of Theorem \ref{thm}}
%%%%%%%%%%%%%%%%%%%%%%%%%%%%%%%%%%%%%%%%%%%%%%%%%%%%
%%%%%%%%%%%%%%%%%%%%%%%%%%%%%%%%%%%%%%%%%%%%%%%%%%%%

We will use the prime number theorem in the proof of Theorem \ref{thm}.

\begin{pnt}[\cite{Dittrich},\cite{RJ}]
Let $\pi(x)$ be the number of primes less than or equal to $x$. Then
$$\lim_{x\to\infty}\frac{\pi(x)}{x/\ln x}=1.$$
\end{pnt}

\begin{lem}\label{lem}
Recall that $y>0$. Let $0\leq\gamma<2\pi$.
There are at most two primes $p$ such that
$$y\ln p=2n\pi+\gamma$$
for some $n\in\mathbb{N}\cup\{0\}.$
\end{lem}

\begin{proof}
Suppose that there exist three distinct primes $p_1<p_2<p_3$ and $\ell,m,n\in\mathbb{N}\cup\{0\}$ such that
\begin{equation}\label{eq1}
y\ln p_1=2\ell\pi+\gamma,\quad y\ln p_2=2m\pi+\gamma,\quad y\ln p_3=2n\pi+\gamma.
\end{equation}
This leads to a contradiction. From eq. (\ref{eq1}), we have
\begin{equation}\label{eq2}
y(\ln p_2-\ln p_1)=2(m-\ell)\pi,\qquad  y(\ln p_3-\ln p_1)=2(n-\ell)\pi.
\end{equation}
Notice that $\ell<m<n$. Let $m-\ell=h$ and $n-\ell=k$. From eq. (\ref{eq2}), we have
$$\frac{\ln p_3-\ln p_1}{\ln p_2-\ln p_1}=\frac{k}{h}.$$
Therefore
$$h(\ln p_3-\ln p_1)=k(\ln p_2-\ln p_1)$$
and hence
$$\left(\frac{p_3}{p_1}\right)^h=\left(\frac{p_2}{p_1}\right)^k.$$
Thus
$$p_1^kp_3^h=p_1^hp_2^k.$$
This contradicts the uniqueness of prime factorization.
\end{proof}

\begin{defn}\label{defn}
Recall that $y>0$ and $0<K<1$. Let $\beta$ be the number such that
$$\cos \beta=K,\quad 0<\beta<\frac{\pi}{2}.$$
For each $n\in\mathbb{N}\cup\{0\}$, let
\begin{eqnarray*}
A_n&=&\left\{p\in P\mid 2n\pi-\beta<y\ln p+\alpha\leq 2n\pi+\beta\right\},\\
B_n&=&\left\{p\in P\mid (2n+1)\pi-\beta<y\ln p+\alpha\leq (2n+1)\pi+\beta\right\}
\end{eqnarray*}
and
$$A=\bigcup_{n=0}^\infty A_n,\qquad B=\bigcup_{n=0}^\infty B_n.$$
\end{defn}

\noindent
{\textbf{\large{Proof of Theorem \ref{thm}}}}\\

\noindent
Notice that $P^+\subset A$ and $P^-\subset B$. From Lemma \ref{lem}, we know that
$A-P^+$ has at most two elements and $B-P^-$ also has at most two elements. Therefore, it is enough to show that
$$\sum_{p\in A}\frac{1}{p}=\infty\quad{and}\quad  \sum_{p\in B}\frac{1}{p}=\infty.$$

\indent
Recall that $y>0$. By the prime number theorem, there exists $M>0$ such that if $x>M$ then
\begin{equation}\label{eq:good}
e^{-\frac{\beta}{2y}}\frac{x}{\ln x}\leq\pi(x)\leq  e^{\frac{\beta}{2y}}\frac{x}{\ln x}.
\end{equation}
From Definition \ref{defn}, we have
$$A_n=\left\{p\in P\mid e^{\frac{2n\pi}{y}-\frac{\beta+\alpha}{y}}< p\leq e^{\frac{2n\pi}{y}+\frac{\beta-\alpha}{y}}\right\}$$
and
$$B_n=\left\{p\in P\mid e^{\frac{(2n+1)\pi}{y}-\frac{\beta+\alpha}{y}}< p\leq e^{\frac{(2n+1)\pi}{y}+\frac{\beta-\alpha}{y}}\right\}.$$

Notice that $A_1,B_1,A_2,B_2,\ldots$ are mutually disjoint. There exists $N\in\mathbb{N}$ such that if $n>N$ then
$$e^{\frac{2n\pi}{y}-\frac{\beta+\alpha}{y}}>M.$$
From now on, we assume that $n>N$. By eq. (\ref{eq:good}), we can find lower bounds for the number of elements in $A_n$ and $B_n$.
We have
\begin{eqnarray}\label{eq:An}
|A_n|&\geq& e^{-\frac{\beta}{2y}}\frac{e^{\frac{2n\pi}{y}+\frac{\beta-\alpha}{y}}}{\frac{2n\pi}{y}+\frac{\beta-\alpha}{y}}-
e^{\frac{\beta}{2y}}\frac{e^{\frac{2n\pi}{y}-\frac{\beta+\alpha}{y}}}{\frac{2n\pi}{y}-\frac{\beta+\alpha}{y}}  \nonumber\\
&&\qquad\qquad\quad
=\frac{ye^{\frac{2n\pi}{y}+\frac{\beta-2\alpha}{2y}}}{2n\pi+\beta-\alpha}-\frac{ye^{\frac{2n\pi}{y}-\frac{\beta+2\alpha}{2y}}}{2n\pi-\beta-\alpha}
\end{eqnarray}
and
\begin{eqnarray}\label{eq:Bn}
|B_n|&\geq& e^{-\frac{\beta}{2y}}\frac{e^{\frac{(2n+1)\pi}{y}+\frac{\beta-\alpha}{y}}}{\frac{(2n+1)\pi}{y}+\frac{\beta-\alpha}{y}}-
e^{\frac{\beta}{2y}}\frac{e^{\frac{(2n+1)\pi}{y}-\frac{\beta+\alpha}{y}}}{\frac{(2n+1)\pi}{y}-\frac{\beta+\alpha}{y}} \nonumber\\
&&\qquad\qquad\quad =\frac{ye^{\frac{(2n+1)\pi}{y}+\frac{\beta-2\alpha}{2y}}}{(2n+1)\pi+\beta-\alpha}
-\frac{ye^{\frac{(2n+1)\pi}{y}-\frac{\beta+2\alpha}{2y}}}{(2n+1)\pi-\beta-\alpha}.
\end{eqnarray}

\indent
Notice that, if $p\in A_n$, then
\begin{equation}\label{eq:AnR}
\frac{1}{p}\geq e^{-\frac{2n\pi}{y}-\frac{\beta-\alpha}{y}}
\end{equation}
and if $p\in B_n$ then
\begin{equation}\label{eq:BnR}
\frac{1}{p}\geq e^{-\frac{(2n+1)\pi}{y}-\frac{\beta-\alpha}{y}}.
\end{equation}

From eq. (\ref{eq:An}) and (\ref{eq:AnR}), we have
\begin{eqnarray*}
\sum_{p\in A_n}\frac{1}{p}&\geq& \left(\frac{ye^{\frac{2n\pi}{y}+\frac{\beta-2\alpha}{2y}}}{2n\pi+\beta-\alpha}
-\frac{ye^{\frac{2n\pi}{y}-\frac{\beta+2\alpha}{2y}}}{2n\pi-\beta-\alpha}\right)
 e^{-\frac{2n\pi}{y}-\frac{\beta-\alpha}{y}}\\
&=&\frac{ye^{-\frac{\beta}{2y}}}{2n\pi+\beta-\alpha}-\frac{ye^{-\frac{3\beta}{2y}}}{2n\pi-\beta-\alpha}\\
&=&y\frac{(2n\pi-\beta-\alpha)e^{-\frac{\beta}{2y}}-(2n\pi+\beta-\alpha)e^{-\frac{3\beta}{2y}}}{(2n\pi-\alpha)^2-\beta^2}\\
&=&y\frac{(2n\pi-\alpha)\left(e^{-\frac{\beta}{2y}}-e^{-\frac{3\beta}{2y}}\right)-\beta\left(e^{-\frac{\beta}{2y}}+e^{-\frac{3\beta}{2y}}\right)}
{(2n\pi-\alpha)^2-\beta^2}\\
&=&\frac{2cn-d}{(2n\pi-\alpha)^2-\beta^2}.
\end{eqnarray*}
where
\begin{equation}\label{eq:c}
c=y\pi\left(e^{-\frac{\beta}{2y}}-e^{-\frac{3\beta}{2y}}\right)>0
\end{equation}
and
$$d=y\alpha\left(e^{-\frac{\beta}{2y}}-e^{-\frac{3\beta}{2y}}\right)+y\beta \left(e^{-\frac{\beta}{2y}}+e^{-\frac{3\beta}{2y}}\right).$$

\noindent
Similarly from eq. (\ref{eq:Bn}) and (\ref{eq:BnR}), we have
\begin{eqnarray*}
\sum_{p\in B_n}\frac{1}{p}&\geq& \left(\frac{ye^{\frac{(2n+1)\pi}{y}+\frac{\beta-2\alpha}{2y}}}{(2n+1)\pi+\beta-\alpha}
-\frac{ye^{\frac{(2n+1)\pi}{y}-\frac{\beta+2\alpha}{2y}}}{(2n+1)\pi-\beta-\alpha}\right)
 e^{-\frac{(2n+1)\pi}{y}-\frac{\beta-\alpha}{y}}\\
&=&\frac{ye^{-\frac{\beta}{2y}}}{(2n+1)\pi+\beta-\alpha}-\frac{ye^{-\frac{3\beta}{2y}}}{(2n+1)\pi-\beta-\alpha}\\
&=&y\frac{((2n+1)\pi-\beta-\alpha)e^{-\frac{\beta}{2y}}-((2n+1)\pi+\beta-\alpha)e^{-\frac{3\beta}{2y}}}{(2n\pi+\pi-\alpha)^2-\beta^2}\\
&=&y\frac{((2n+1)\pi-\alpha)\left(e^{-\frac{\beta}{2y}}-e^{-\frac{3\beta}{2y}}\right)-\beta\left(e^{-\frac{\beta}{2y}}+e^{-\frac{3\beta}{2y}}\right)}
{(2n\pi+\pi-\alpha)^2-\beta^2}\\
&=&\frac{c(2n+1)-d}{(2n\pi+\pi-\alpha)^2-\beta^2}.
\end{eqnarray*}

Recall eq. (\ref{eq:c}). Since $c>0$, we have
$$\sum_{p\in A}\frac{1}{p}\geq\sum_{n=N+1}^\infty\sum_{p\in A_n}\frac{1}{p}\geq \sum_{n=N+1}^\infty \frac{2cn-d}{(2n\pi-\alpha)^2-\beta^2}=\infty$$
and
$$\sum_{p\in B}\frac{1}{p}\geq\sum_{n=N+1}^\infty\sum_{p\in B_n}\frac{1}{p}\geq \sum_{n=N+1}^\infty \frac{c(2n+1)-d}{(2n\pi+\pi-\alpha)^2-\beta^2}=\infty.$$
Thus
$$\sum_{p\in A}\frac{1}{p}=\infty\quad{and}\quad  \sum_{p\in B}\frac{1}{p}=\infty.$$
\qed

%%%%%%%%%%%%%%%%%%%%%%%%%%%
%%%%%%%%%%%%%%%%%%%%%%%%%%%

\end{document}